\documentclass[11pt,a4paper]{amsart}
\usepackage{amsmath,amsfonts,amsthm,amsopn,color,amssymb,enumitem}
\usepackage[a4paper]{geometry}
\usepackage{palatino}
\usepackage{graphicx}
\usepackage[colorlinks=true]{hyperref}
\hypersetup{urlcolor=blue, citecolor=red, linkcolor=blue}

\usepackage{cite}
\usepackage{relsize}
\usepackage{esint}
\usepackage{verbatim}
\usepackage{mathrsfs}
\usepackage{xcolor}
\usepackage{tikz}

\newcommand{\e}{\varepsilon}

\newcommand{\R}{\mathbb{R}}

\newcommand{\RN}{{\mathbb{R}^N}}

\DeclareMathOperator{\dist}{dist}

\renewcommand{\leq}{\leqslant}
\renewcommand{\geq}{\geqslant}

\newcommand{\bd}{\boldsymbol}

\newcommand{\beq}{\begin{equation}}
\newcommand{\eeq}{\end{equation}}

\newtheorem{theorem}{Theorem}[section]
\newtheorem*{theorem*}{Theorem}
\newtheorem{lemma}[theorem]{Lemma}

\newtheorem{definition}[theorem]{Definition}
\newtheorem{proposition}[theorem]{Proposition}
\newtheorem{corollary}[theorem]{Corollary}
\theoremstyle{definition}
\newtheorem{remark}[theorem]{Remark}

\begin{document}

\title[Sign-changing bubbling solutions on polygons]{Sign-changing bubbling solutions on polygons for the Brezis-Nirenberg problem in high dimensions}

\author[G. M. Rago]{Giuseppe Mario Rago}
\address{\noindent  Dipartimento di Matematica, Universitá degli Studi di Bari Aldo Moro, Italy }
\email{g.rago6@phd.uniba.it}

\subjclass{35B44, 35B33, 35J25}
\keywords{Brezis-Nirenberg problem; blowing-up solutions; sign-changing solutions; Lyapunov-Schmidt reduction.}

\begin{abstract}
We consider the Brezis--Nirenberg problem $-\Delta u = |u|^{2^{*}-2}u + \e u$ in a smooth bounded domain $\Omega \subset \mathbb{R}^{N}$, $N \geq 5$, with Dirichlet boundary conditions, where $\e > 0$ is small. Two multi-peak configurations found in dimension four in \cite{PRV2} persist. In an annulus, for every even $k \geq 2$, we construct solutions blowing up at the vertices of a regular $k$-gon with alternating signs. In the ball, we construct solutions blowing up positively at the vertices of a regular $k$-gon in a coordinate plane and negatively at those of a regular $k$-gon in an orthogonal plane, for $2 \leq k \leq 5$ if $N = 5,6$ and $2 \leq k \leq 4$ if $N \geq 7$. This range is sharp: for larger $k$ no such solution exists with parameters in a compact set, while otherwise the reduced energy has a unique, non-degenerate critical point.
\end{abstract}

\maketitle

\section{Introduction and main results}

In this paper we continue the study, started in \cite{PRV2}, of sign-changing blowing-up solutions of the Brezis--Nirenberg problem
\begin{equation}\label{pb}
  \begin{cases}
    -\Delta u=|u|^{2^*-2}u+\varepsilon u & \text{in }\Omega,\\
    u=0 & \text{on }\partial\Omega,
  \end{cases}
\end{equation}
where $\Omega\subset\mathbb R^N$ with $N\geq 5$ is a smooth bounded domain, $\varepsilon>0$ is a small parameter and, as usual, $2^*=\frac{2N}{N-2}$ is the critical exponent for the Sobolev embedding. \\
Problem \eqref{pb} was introduced by Brezis and Nirenberg in their celebrated paper \cite{BN}.\\
A particular feature of \eqref{pb}, due to the critical behaviour of the nonlinearity, is the possible existence of solutions which blow up at one or more points in the domain as the parameter $\e$ approaches $0$ in dimension $N\geq 4$ or the parameter $\e$ approaches a positive number $\lambda_*$ in dimension $N=3$. \ \\
The profile of the positive blowing-up solutions and the construction of such positive solutions predicted by the asymptotic analysis have been the subject of a wide
literature (see for instance \cite{D}, \cite{H}, \cite{KL1}, \cite{KL2}, \cite{MRV}, \cite{MP}, \cite{PRV}, \cite{R}  and references therein). \\ Instead, the profile of the sign-changing solutions is far from being completely understood. As far as we know the complete scenario of sign-changing solutions has been obtained only in
the radial case by Esposito, Ghoussoub, Pistoia and Vaira in \cite{EGPV} in dimensions $N\geq 7$
and by Amadori, Gladiali, Grossi, Pistoia and Vaira \cite{AGGPV} in lower dimensions.\\
Sign-changing solutions for \eqref{pb} have been built by Micheletti and Pistoia \cite{MiPi}, Iacopetti and Vaira  \cite{IV1, IV2}, Li, Vaira, Wei and Wu \cite{LVWW, LVWW1}, Pistoia and Vaira in \cite{PV} and more recently in \cite{PRV2} and \cite{PRV3}.
\\\\
Indeed, in \cite{PRV2} Pistoia, Rago and Vaira constructed sign-changing solutions with several concentration points in dimension $N=4$, showing that the locations and the relative scales of the bubbles are governed by a signed Green–Robin interaction matrix. They formulated an abstract existence criterion in terms of a simple positive eigenvalue with a positive eigenvector and of a stable critical set and applied it to a positive-negative pair in a general domain and to several symmetric multi-peak configurations: alternating regular polygons, orthogonal polygons, one central peak surrounded by peaks of the opposite sign and aligned patterns with three, four and five peaks. For the two-peak solution they also proved that it has exactly two nodal domains and, under a natural balance condition and the connectedness of the boundary, that the closure of its nodal set meets the boundary. \ \\
Instead, in \cite{PRV3} the authors considered the configuration with a positive peak at the centre of the unit ball surrounded by $k$ negative peaks at the vertices of a regular polygon. In dimension $N=5$ such solutions exist for every $k \geq 2$ and the polygon approaches the boundary of the ball as $k\to+\infty$ at distance $\mathcal{O}(\frac{1}{k})$. \ \\
On the other hand, for $N\geq 6$ there are two distinct families whenever $2\leq k \leq6$, whereas for $k$ large no solution of this type exists with parameters in a compact set. Moreover, for $2\leq k \leq6$ the nodal set of these solutions always meets the boundary of the ball. \ \\
In dimension $N\geq 5$, sign-changing multi-bubble solutions for critical and supercritical problems were studied by Micheletti and Pistoia in \cite{MiPi}. In particular, for the Brezis-Nirenberg perturbation in dimensions $N\ge5$, their reduction leads to a finite-dimensional functional depending simultaneously on the
concentration points and on the individual scales.
In order to understand the difficulties and also the main results of this paper, we introduce some notation used throughout the paper. \\
Let $G_\Omega(x, y)$ be the Green function of $-\Delta$ in $\Omega$ with Dirichlet boundary conditions and let $H$ be its regular part, i.e. 
 $$G_\Omega(x, y):=\frac{1}{\omega_N}\frac{1}{|x-y|^{N-2}}-H(x, y)$$ where $\omega_N:=(N-2){\rm meas}(\mathbb S^{N-1})$, $\mathbb S^{N-1}$ is the unit sphere of $\mathbb{R}^{N}$ and for every $y\in\Omega$, $$-\Delta_{x} H(x, y)=0\ \hbox{in}\ \Omega,\ H(x, y)=\frac{1}{\omega_N|x-y|^{N-2}}\ \hbox{on}\  \partial\Omega.$$  
The leading term of the regular part, $\tau_\Omega(x)=H(x,x)$, is called the \emph{Robin function} of $\Omega$. \ \\
Given an integer $k\geq 2$ we set
$$\Omega^*_k:=\left\{\boldsymbol\xi:=(\xi_1, \ldots, \xi_k)\in \Omega^k\,\,:\,\, \xi_i\neq \xi_j \,\, \mbox{for}\,\, i\neq j\right\}$$
and, for any $\boldsymbol\xi\in\Omega^*_k$ and $\boldsymbol\sigma=(\sigma_1,\ldots,\sigma_k)\in\{-1,1\}^k$, we define the symmetric matrix 
$$
M^{\boldsymbol\sigma}(\boldsymbol\xi)=
\begin{pmatrix}
\tau_\Omega(\xi_1)&-\sigma_1\sigma_2G_{\Omega}(\xi_1,\xi_2)&\cdots&-\sigma_1\sigma_kG_{\Omega}(\xi_1,\xi_k)\\
-\sigma_1\sigma_2G_{\Omega}(\xi_1,\xi_2)&\tau_\Omega(\xi_2)&\cdots&-\sigma_2\sigma_kG_{\Omega}(\xi_2,\xi_k)\\
\vdots&\vdots&\ddots&\vdots\\
-\sigma_1\sigma_kG_{\Omega}(\xi_1,\xi_k)&-\sigma_2\sigma_kG_{\Omega}(\xi_2,\xi_k)&\cdots&\tau_\Omega(\xi_k)
\end{pmatrix}.$$
The matrix $M^{\boldsymbol\sigma}(\boldsymbol\xi)$ is the Green--Robin interaction matrix: its diagonal entries describe the interaction of each bubble with the boundary, while the off-diagonal terms describe the interaction between distinct bubbles. \ \\
In our notation, the leading part of the reduced energy in \cite{MiPi} is
\begin{equation}\label{signed-reduced-functional}
F^{\boldsymbol\sigma}_k(\boldsymbol d,\boldsymbol\xi)= \frac{a_N}{2}\left[\sum_{j=1}^k \tau_\Omega(\xi_j) d_j^{N-2}-\sum_{i, j=1\atop  i\neq j}^k \sigma_j\sigma_i G_{\Omega}(\xi_i, \xi_j)d_i^{\frac{N-2}{2}}d_j^{\frac{N-2}{2}}\right]-b_N\sum_{j=1}^k d_j^2\end{equation}
or equivalently, setting $\kappa_j:=d_j^{\frac{N-2}{2}}$,
\begin{equation}\label{signed-reduced-functional1}
F^{\boldsymbol\sigma}_k(\boldsymbol \kappa,\boldsymbol\xi)=\frac{a_{N}}{2} \left\langle
M^{\boldsymbol\sigma}(\boldsymbol\xi)\boldsymbol \kappa,
\boldsymbol \kappa
\right\rangle
-b_N\sum_{j=1}^k\kappa_j^{\frac{4}{N-2}},
\end{equation}
where $a_N, b_N>0$ and explicitly given in \eqref{aNbN}.\\
Thus the geometry of the domain and the nodal pattern enter the reduced
problem simultaneously.\\
The general reduction allows one to prescribe arbitrary signs, but the critical point analysis of \eqref{signed-reduced-functional} is rather delicate because the location variables and all the relative scales remain coupled. This is one of the reasons why explicit constructions in natural geometries were obtained only for particular configurations. For example, in the ball, in \cite{MiPi} Micheletti and Pistoia constructed a family with an even number of bubbles arranged on a regular polygon with alternating signs; the symmetry reduces all the scales to one parameter and all the locations to one radial variable. The same paper also produces many sign-changing patterns in specially constructed dumbbell-type domains and also the case of one positive and one negative bubble in an arbitrary smooth bounded domain is considered. Here the difficulty is not the general reduction itself, but the analysis of its finite-dimensional functional and the associated critical points. In dimension $N=4$ the situation is different: the concentration scales are exponentially small in $1/\varepsilon$ and the leading equations for the relative scales become spectral.  Hence the scale variables can be eliminated through the spectral structure of the interaction matrix which makes accessible configurations that are considerably harder to detect through the higher-dimensional reduced functional. 
In this paper we show that two of the multi-peak configurations constructed in \cite{PRV2} survive in dimensions $N\ge5$ : alternating regular polygons in annuli and two regular polygons lying in orthogonal coordinate planes. \\ \\
For $r>0$ and $k\geqslant2$ we set $\theta_j:=\frac{2\pi(j-1)}{k}$ and
\begin{equation}\label{points}
\begin{aligned}
& \xi_j(r):=r(\cos\theta_j,\sin\theta_j,0,\dots,0), \\&
\eta_j(r):=r(0,0,\cos\theta_j,\sin\theta_j,0,\dots,0),\qquad j=1,\dots,k.
\end{aligned}
\end{equation} \ \\
Our main results are the following. \ \\ 
\begin{theorem}[Symmetric multi-peak configurations]
\label{main} Let $N\geq 5$.
For $\varepsilon>0$ sufficiently small, problem~\eqref{pb} has the
following sign-changing solutions.

\begin{enumerate}
\item[(a)] \textbf{Alternating regular polygons in an annulus.}
Let 
\begin{equation}\label{omegarho}
\Omega_\rho=\{x\in\RN:\rho<|x|<1\}
\end{equation}
with $0<\rho<1$. For every even $k\geq2$ there is a solution
$$u_{\e} = \sum_{j=1}^{k} (-1)^{j+1} P\mathcal U_{\delta_{\e}, \xi_{j}(r_{\e})} + \phi_{\e}, \ \ \mbox{with} \ \delta_{\e}=d_{\e}\e^{\frac{1}{N-4}} \ \ \mbox{and} \ \ \ \| \phi_{\e} \|_{H^{1}_{0}(\Omega_{\rho})} \to 0$$
 with $k$ concentration points of alternating sign at the vertices of a regular $k$-gon in the
$(x_1,x_2)$-plane. It is even in $x_2$ and in $x_3,\dots,x_N$, and anti-invariant under the
rotation $R_{2\pi/k}$ of angle $2\pi/k$ in the $(x_1,x_2)$-plane: $u(R_{2\pi/k}x)=-u(x)$. Moreover, up to a subsequence, $r_{\e} \to r_{0}$ and $d_{\e} \to d_0>0$, where $r_0 \in (\rho,1)$ is a global minimum point of the function $\Psi_{k}$ defined in \eqref{Psianello}. \ \\
\item[(b)] \textbf{Two orthogonal regular polygons in the ball.}
Let $\Omega=\mathbb{B}_1(0)$ be the unit ball and let $2\leq k\leq k^*(N)$, where 
\[
k^*(N):=\begin{cases}5 & \text{if } N=5,6,\\ 4 & \text{if } N\geq7.\end{cases}
\]
There is a solution 
$$u_{\e} = \sum_{j=1}^{k} P\mathcal U_{\delta_{\e}, \xi_{j}(r_{\e})} -  \sum_{j=1}^{k} P\mathcal U_{\delta_{\e}, \eta_{j}(r_{\e})} + \phi_{\e}, \ \ \mbox{with} \ \delta_{\e}=d_{\e}\e^{\frac{1}{N-4}} \ \ \mbox{and} \ \ \ \| \phi_{\e} \|_{H^{1}_{0}(\mathbb{B}_{1}(0))} \to 0$$
with $k$ positive concentration points at the vertices of a regular $k$-gon in the $(x_1,x_2)$-plane and $k$ negative concentration points at the vertices of a regular $k$-gon in the $(x_3,x_4)$-plane. It is invariant under the rotations of angle $2\pi/k$ in each of the two planes, even with respect to $(x_2, x_4) \to (-x_2, -x_4)$, radial in $(x_5,\dots,x_N)$ and antisymmetric under the exchange of the two planes:
\[u(\omega,y,z)=-u(y,\omega,z),\qquad y,\omega \in\R^2,\ z\in\R^{N-4}.\]
Moreover, $r_{\e} \to r_0$ and $d_{\e} \to d_0>0$, where $r_0 \in (0,1)$ is the unique global minimum point of the function $\Phi_{k,N}$ defined in \eqref{eq:Phi-def} (see Theorem \ref{main2}).
\end{enumerate}
\end{theorem}
\begin{theorem}[Sharpness of $k^*(N)$]\label{main2}
Let $N\geqslant5$ and $\Omega=\mathbb{B}_1(0)$.
\begin{enumerate}
\item[(i)] If $2\leqslant k\leqslant k^*(N)$, the limiting reduced  functional, restricted to the two-parameter symmetric family, has exactly one critical
point, which is non-degenerate; in particular the solutions of Theorem~\ref{main}\,(b)
satisfy $r_\varepsilon\to r_0$ and $d_\varepsilon\to d_0$ as $\varepsilon\to0$, where $r_0$ is the
unique critical point of $\Phi_{k,N}$.
\item[(ii)] If $k>k^*(N)$, for every compact set $\mathcal K\subset(0,1)\times(0,\infty)$
there is $\varepsilon_{\mathcal K}>0$ such that, for $\varepsilon\in(0,\varepsilon_{\mathcal K})$,
problem \eqref{pb} has no solution of the form given in Theorem~\ref{main}\,(b) with
$(r,d)\in\mathcal K$ (see Theorem~\ref{thm:nonex}).
\end{enumerate}
\end{theorem}
\noindent
The limit scale $d_0$ is given explicitly in Proposition \ref{propabstract} in terms of the minimum value of $\Psi_k$, respectively $\Phi_{k,N}$. 
Moreover, Theorem \ref{main} $(a)$ extends to annuli the alternating polygons constructed in the ball in \cite{MiPi} (see Theorem \ref{TeoMi2}) and we point out that no restriction on the (even) number of peaks is needed. \ \\
Instead, in Theorem \ref{main} $(b)$, the number of peaks is bounded: the restriction $k\leq k^{*}(N)$ comes from the sign of the constant $\mathfrak{c}_{k,N}$ in \eqref{eq:ckN}, which measures the competition between the mutual attraction of the peaks of the same polygon and the repulsion between the two polygons (see Proposition \ref{prop:sign}) and it is optimal by Theorem \ref{main2}. In dimension four the analogous condition gives $k \leq 6$ (see Remark \ref{rmk4D}).

%\subsection{Open problems and further directions.} The results above suggest some natural questions.
%\begin{itemize}
%\item [1.] \emph{Orthogonal polygons with many peaks.} For $k>k_{*}(N)$ the function $\Phi_{k,N}$ tends to $-\infty$ as $r\to0^{+}$ and, for the values of $(k,N)$ we checked numerically, it is strictly increasing on $(0,1)$ (See Remark \ref{rem:negative}). A proof of this monotonicity would give, as in \cite{RV}, (Theorem 1.2 (ii)), the non-existence of solutions of this form with parameters in compact sets. \ \\
%\item [2.] \emph{Uniqueness and non-degeneracy.} By Lemma \ref{lem:critPsi}, a non-degenerate global minimum point of $\Psi_{k}$, respectively $\Phi_{k,N}$, gives a non-degenerate critical point of the reduced energy. It would be interesting to prove that the global minimum is unique and non-degenerate, which would give the local uniqueness of the solutions. \ \\
%\item [3.] \emph{The nodal set.} As in \cite{PRV2} and \cite{RV}, one may ask whether the nodal set of the solutions of Theorem \ref{main} meets the boundary of the domain.
%\end{itemize}

\subsection{Organization of the paper}
In Section \ref{2} we recall the Lyapunov–Schmidt reduction and we perform it in a class of symmetric functions; we show that for symmetric configurations the reduced energy is governed by a single function of the radius and we prove an abstract existence result (Proposition \ref{propabstract}). In Section \ref{anello} we analyze the case of the annulus, proving Theorem \ref{main} (a) and finally, in Section \ref{palla} we study the case of the two orthogonal regular polygons in the ball, proving Theorem \ref{main} (b) and Theorem \ref{main2}.
\ \\ \\
\noindent{\bf Acknowledgments.}
The author is deeply grateful to Professor Angela Pistoia and Professor Giusi Vaira for many stimulating discussions, valuable insights and helpful suggestions,
which greatly contributed to the development of this work. \ \\
Moreover the author acknowledges support from the INdAM-GNAMPA project ``Singolarità, nonlinearità e interazioni non locali: nuovi approcci analitici per equazioni differenziali e sistemi'' (CUP E53C25002010001). \ \\
The author acknowledges also the use of ChatGPT (OpenAI) to assist selected calculations. He has independently verified all mathematical arguments and calculations to take full responsibility for the content of the paper.

\section{Preliminaries}\label{2}
\noindent First we recall the definition of  a stable critical point we use in this paper (see \cite{L}).
\begin{definition} \label{stable}
Let $f:D \to \mathbb R$ be a $C^1$-function, where $D \subset \mathbb R^N$ is an open set.\\ We say that
$\xi_0$ is a stable critical point of $f$ if $\nabla f(\xi_0)=0$ and there exists a neighbourhood $V\subset\subset D$ of $\xi_0$ such that $\nabla f(\xi)\neq 0$ for any $\xi\in\partial V$,
$$\nabla f(\xi)=0,\quad \xi\in V \Longrightarrow f(\xi)=f(\xi_0)$$
and
$${\rm deg }(\nabla f,V , 0)\neq  0,$$ where deg denotes the Brouwer degree.
\end{definition}
\noindent Next, let us define the {\it bubble} (see \cite{A, CGS, T})
\beq\label{bubble}
\mathcal U_{\delta, \xi}(x)=\alpha_N \frac{\delta^{\frac{N-2}{2}}}{\left(\delta^2+|x-\xi|^2\right)^{\frac{N-2}{2}}},\quad x, \xi\in\mathbb R^N,\,\, \delta>0\eeq where $\alpha_N=\left[N(N-2)\right]^{\frac{N-2}{4}}$. As is well-known, these are all the positive solutions of the problem \beq\label{critico}-\Delta \mathcal U=\mathcal U^{\frac{N+2}{N-2}},\quad\mbox{in}\,\, \mathbb R^N.\eeq 
Moreover, we denote by $P\mathcal U_{\delta, \xi}$ the projection of the bubble into $H^1_0(\Omega)$, i.e. 
$$\Delta P\mathcal U_{\delta, \xi}=\Delta \mathcal U_{\delta, \xi}\quad \mbox{in}\,\Omega, \qquad P\mathcal U_{\delta, \xi}=0\quad\mbox{on}\,\, \partial\Omega$$
and we set $\mathcal U : = \mathcal U_{1,0}$ and
%Now, let $k\geq 2$ and let $\bd\xi\in\Omega^*_k$. Let $\eta\in (0, 1)$ and set $$\mathcal O_\eta:=\left\{(\bd d, \bd\xi)\in\mathbb R^k_+\times \Omega_k^*\,\, :\,\,d_j\in (\eta, \eta^{-1}),\,\, {\rm dist}(\xi_j, \partial\Omega)\geq \eta,\,\, |\xi_j-\xi_h|\geq \eta,\,\, j\neq h\right\}.$$
%Let us also define 
\beq\label{aNbN}
a_N:=\int_{\mathbb R^N}\mathcal U^{2^*-1}(x)\, dx,\qquad b_N:=\frac 12 \int_{\mathbb R^N}\mathcal U^2(x)\, dx\eeq
which are finite since $N \geq 5$. \ \\ \\
Moreover, in \cite{MiPi}, the following result is proved.

\begin{theorem}\label{TeoMiPi}(\cite[Theorem 3.4]{MiPi})
Let $(\bd d_0, \bd\xi_0)\in (0,+\infty)^{k} \times \Omega_{k}^{*}$ be a stable critical point of the function $F^{\bf\sigma}_k(\bd d, \bd \xi)$ defined in \eqref{signed-reduced-functional} in the sense of Definition \ref{stable}.\\
Then for $\e$ small enough, problem \eqref{pb} has a family of solutions 
$$u_\e = \sum_{j=1}^k \sigma_j P \mathcal U_{\delta_j^\e, \xi_j^\e}+\phi_{\e}, \ \ \ \ \ \ \ \ \| \phi_{\e} \|_{H^{1}_{0}(\Omega)} \to 0, $$
such that $\delta_j^\e=d_j^\e \e^{\frac{1}{N-4}}$; moreover, up to a subsequence $d_j^\e \to d_{j}^*$ and $\xi_j^\e\to\xi_{j}^*$ where $(\bd d^*, \bd \xi^*)$ is a critical point of  $F^{\bf\sigma}_k$ with  $F^{\bf\sigma}_k(\bd d_0, \bd \xi_0)= F^{\bf\sigma}_k(\bd d^*, \bd\xi^*)$. 
\end{theorem}

\noindent Theorem \ref{TeoMiPi} establishes a general criterion to find a sign-changing solution of problem \eqref{pb} in dimension $N\geq 5$. It remains to find domains and sign patterns for which the function $F^{\bf\sigma}_k$ has a stable critical point. In \cite{MiPi} this is done in two situations. The first one is a dumbbell-shaped domain with thin handles, in which problem \eqref{pb} has many sign-changing solutions with many nodal domains.
\begin{theorem}\label{TeoMi1}(\cite{MiPi}, Theorem 1.4)
Assume $N\geq 5$ and let $k\geq 2$ be a fixed integer. There exists a contractible domain $\Omega$
such that for $\e$ small enough, problem \eqref{pb} has at least $2^{k-1}-1$ pairs of sign-changing
solutions that blow up positively at $h$  points and blow up negatively at $k-h$ points for any
$h= 1, \ldots, k-1$.\end{theorem}
\noindent The second one is the ball, with peaks of alternating signs.
\begin{theorem}\label{TeoMi2}(\cite{MiPi}, Theorem 1.5)
Assume $\Omega$ is the unit ball in $\mathbb R^N$ with $N\geq 5$. Let $h \geq 1$. Then for $\e$ small enough,
problem \eqref{pb} has at least a pair of sign-changing solutions that blow up positively at $h$ points
and blow up negatively at $h$ points.\end{theorem}
\noindent In both cases the symmetries are used to reduce the problem to a function in which the location and the scale variables are not fully coupled. Our aim is to find other configurations of sign-changing solutions, extending some of the results proved in dimension $N = 4$ in \cite{PRV2}. \ \\
We recall that the solutions of \eqref{pb} are the critical points of the following functional
$$J_{\e}(u)=\frac 12 \int_{\Omega} |\nabla u|^{2} \ dx -\frac{\e}{2} \int_{\Omega} u^{2} \ dx - \frac{1}{2^{*}} \int_{\Omega} |u|^{2^{*}} \ dx, \quad u \in H^{1}_{0}(\Omega).$$

\subsection{The symmetric setting}\label{subsec:sym}
Theorem~\ref{TeoMiPi} requires a stable critical point of $F^\sigma_k$ in the whole space $(0,+\infty)^k\times\Omega^*_k$, whereas the configurations of Theorem~\ref{main} are critical points of $F^\sigma_k$ restricted to a two-dimensional set of symmetric parameters and the stability, in the sense of Definition~\ref{stable}, is not inherited by the restriction. For this reason we perform the reduction directly in a space of symmetric functions, in which only two parameters survive. We treat the two parts of Theorem~\ref{main} at the same time.

\medskip
\noindent\textbf{The two configurations.} 
\begin{enumerate}[label=\textup{(\Alph*)},leftmargin=*]
\item \emph{Alternating regular polygons in an annulus.} Let us consider the annulus $\Omega_{\rho}$ as in \eqref{omegarho}, an even integer $k\geqslant2$ and we set $I:=(\rho,1)$ and $\ell:=k$. \ \\
Denoting by $R$ the rotation of angle $\frac{2\pi}k$ in the $(x_1,x_2)$-plane, the space of symmetric functions is
\[\mathcal H_{s}:=\bigl\{u\in H^1_0(\Omega_\rho):\ u(Rx)=-u(x),\ u\text{ is even in }x_h\text{ for }h=2,\dots,N\bigr\}.\]

\item \emph{Two orthogonal regular polygons in the ball.} Let $\Omega=\mathbb{B}_{1}(0)$, $k\geqslant2$, $I:=(0,1)$ and $\ell:=2k$. We write
$$x=(y,w,z)\in\R^2\times\R^2\times\R^{N-4},$$
with $y=(x_1,x_2)$, $w=(x_3,x_4)$, $z=(x_5,\dots,x_N)$. The bubbles are centred at $\xi_{1}(r), \cdots, \xi_{k}(r)$ with a positive sign $+1$ and at $\xi_{k+j}(r) :=\eta_{j}(r)$, $j=1, \cdots, k$, with negative sign $-1$; that is, $\sigma_{j}=1$ for $j \leq k$ and $\sigma_{j}=-1$ for $j>k$. \ \\
We introduce the following orthogonal maps
\begin{equation}\label{eq:maps}
\begin{gathered}
R_1(y,w,z)=\bigl(e^{\frac{2\pi i}k}y,w,z\bigr),\qquad R_2(y,w,z)=\bigl(y,e^{\frac{2\pi i}k}w,z\bigr),\qquad S(y,w,z)=(\bar y,\bar w,z),\\
\widehat Q(y,w,z)=(y,w,Qz),\ \ \mbox{with} \ \ Q \in O(N-4) \ \ \ \mbox{and} \ \ T(y,w,z)=(w,y,z).
\end{gathered}
\end{equation}
Therefore, the space of symmetric functions is
\begin{equation}\label{eq:Hs}
\mathcal H_{s}:=\Bigl\{u\in H^1_0(\mathbb{B}_1(0)):\ u\circ R_1=u\circ R_2=u\circ S=u\circ\widehat Q=u\ \ \forall Q\in O(N-4),\ \ u\circ T=-u\Bigr\}.
\end{equation}
In other words, $u\in \mathcal H_{s}$ if and only if $u$ is invariant under the rotations of angle $\frac{2\pi}k$ in each of the two coordinate planes, even with respect to $(x_2,x_4)\mapsto(-x_2,-x_4)$, radial in $z$, and odd under the exchange of the two planes, $u(w,y,z)=-u(y,w,z)$.
\end{enumerate}
In both cases the maps appearing in the definition of $\mathcal H_{s}$ are orthogonal and leave $\Omega$ invariant; hence, the critical points of $J_\e|_{\mathcal H_{s}}$ are critical points of $J_\e$. \ \\ \\
We look for solutions of the form
\begin{equation}\label{ansatz}
u_\e=V_y+\phi,\qquad V_y:=\sum_{j=1}^\ell\sigma_jP\mathcal U_{\delta,\xi_j(r)},\qquad \delta=d\e^{\frac1{N-4}},
\end{equation}
where $y=(r,d)$ lies in a compact subset of $\mathcal M:=I\times(0,\infty)$ and $\phi\in \mathcal H_{s}$ satisfies suitable orthogonality conditions and it is such that $\|\phi\|_{H^1_0(\Omega)}\to0$ as $\e\to0$. For $\eta>0$ small we set
\[
\mathcal P_\eta:=\{r\in I:\ \dist(r,\partial I)\geqslant\eta\}\times[\eta,\eta^{-1}]\subset\mathcal M.
\]
It is immediate to check that $V_y\in \mathcal H_{s}$. The Lyapunov--Schmidt reduction is by now standard and we do not repeat it here: we refer to \cite{BMP,MiPi,MP} for the details and only recall its outcome. Because of the symmetries, as in \cite{PRV3}, the relevant kernel is two-dimensional:
\[Z:=\mbox{span} \left\{\sum_{j=1}^\ell\sigma_j\,\partial_\delta P\mathcal U_{\delta,\xi_j(r)}\Big|_{\delta=d\e^{\frac{1}{N-4}}},\ \sum_{j=1}^\ell\sigma_j\,\partial_rP\mathcal U_{\delta,\xi_j(r)}\right\}\subset \mathcal H_{s} .
\]
We denote by $Z^\perp$ the orthogonal complement of $Z$ in $\mathcal H_{s}$ and by $\Pi^\perp$ the orthogonal projection on it. Finally, we let $\mu=\e^{\frac{1}{N-4}}$ and
\[
\vartheta(\mu):=\begin{cases}\mu^{5/2}, & N=5,\\ \mu^4|\log\mu|^{2/3}, & N=6,\\ \mu^{\frac{N+2}2}, & N\geqslant7.\end{cases}
\]

\begin{proposition}[\cite{BMP,MP}]\label{prop:reduction}
There are $\e_0,C>0$ such that for $\e\in(0,\e_0)$ and every $y\in\mathcal P_\eta$ there is a unique $\phi_y\in Z^\perp$ with $\|\phi_y\|\leqslant C\vartheta(\mu)$ such that
\[\Pi^\perp\left[V_y+\phi_y-(-\Delta)^{-1}\left(|V_y+\phi_y|^{2^*-2}(V_y+\phi_y)+\e(V_y+\phi_y)\right)\right]=0.\]
Moreover, the map $y\mapsto\phi_y$ is of class $C^1$.
\end{proposition}

\noindent Hence the whole problem is reduced to the study of a finite-dimensional energy, described by the following results.

\begin{lemma}[\cite{BMP,MP}]\label{lem:reduced}
Let $\widetilde J_\e(y):=J_\e(V_y+\phi_y)$. If $y\in\operatorname{int}\mathcal P_\eta$ is a critical point of $\widetilde J_\e$, then $V_y+\phi_y$ is a solution of \eqref{pb}.
\end{lemma}

\begin{proposition}[\cite{BMP,MiPi,MP,R}]\label{prop:expansion}
There are constants $A,B>0$, depending only on $N$, such that, as $\e\to0$,
\[
\widetilde J_\e(y)=\ell\,\frac{S^{N/2}}{N}+\mu^{N-2}\bigl(\ell\,\Gamma_{A,B}(y)+\varpi_\e(y)\bigr),\qquad \varpi_\e\to0\ \text{in }C^1(\mathcal P_\eta),
\]
where $S$ is the best Sobolev constant,
\begin{equation}\label{PsiAB}
\Gamma_{A,B}(r,d)=\frac A2\,\Phi(r)\,d^{N-2}-\frac B2\,d^2
\end{equation}
and
\begin{equation}\label{Phigeneral}
\Phi(r):=\tau_\Omega(\xi_1(r))-\sum_{j=2}^\ell\sigma_j\,G_\Omega(\xi_1(r),\xi_j(r)).
\end{equation}
\end{proposition}
\noindent
Let $c:=(B/A)^{\frac1{N-4}}$. Then $\Gamma_{A,B}(r,d)=Bc^2\,\Gamma(r,d/c)$, where
\begin{equation}\label{Psi}
\Gamma(r,d):=\frac12\bigl(\Phi(r)\,d^{N-2}-d^2\bigr).
\end{equation}
Since this change of variables is a positive dilation in $d$, critical points of $\Gamma_{A,B}$ and of $\Gamma$ correspond. \ \\
For every $r\in I$ with $\Phi(r)>0$, we set
\begin{equation}\label{dr}
d(r):=\left(\frac2{(N-2)\Phi(r)}\right)^{\frac1{N-4}},
\end{equation}
which is the unique positive zero of the function
$$d\mapsto\partial_d \Gamma(r,d)=d\left(\frac{N-2}2\Phi(r)\,d^{N-4}-1\right).$$

\begin{lemma}\label{lem:critPsi}
Let $\Phi\in C^2(I)$. Then $(r_0,d_0)\in\mathcal M$ is a critical point of $\Gamma$ if and only if
\[
\Phi(r_0)>0,\qquad \Phi'(r_0)=0\qquad\text{and}\qquad d_0=d(r_0).
\]
Moreover, at such a point 
$$D^2\Gamma(r_0,d_0)=\operatorname{diag}\left(\frac12\Phi''(r_0)\,d_0^{N-2},\,N-4\right).$$
In particular, $(r_0,d_0)$ is non-degenerate if and only if $\Phi''(r_0)\neq0$.
\end{lemma}

\begin{proof}
Since 
$$\partial_r\Gamma(r,d)=\frac12\Phi'(r)\,d^{N-2},$$
the first part follows from the definition of $d(r)$ as in \eqref{dr}. \ \\
At a critical point, 
$$\partial^2_{rd}\Gamma(r_0, d_0)=\frac{N-2}2\Phi'(r_0)\,d_0^{N-3}=0$$
and, using that
\begin{equation}\label{condminphi}
\frac{N-2}2\Phi(r_0)\,d_0^{N-4}=1,
\end{equation} 
it is easy to check that
$$\partial^2_{dd}\Gamma=\frac{(N-2)(N-3)}2\Phi(r_0)\,d_0^{N-4}-1=N-4.$$
\end{proof}

We can now state the abstract existence result which we shall apply in Sections \ref{anello} and \ref{palla}.

\begin{proposition}\label{propabstract}
Let us be in case \textup{(A)} or \textup{(B)} and assume that $\Phi$ defined in \eqref{Phigeneral} satisfies
\begin{equation}\label{hypPhi}
\Phi(r)>0\ \text{for all }r \in I\qquad\text{and}\qquad \Phi(r)\to+\infty\ \text{ as $r$ tends to the endpoints of }I.
\end{equation}
Then there exists $\e_1>0$ such that, for every $\e\in(0,\e_1)$, problem \eqref{pb} has a sign-changing solution 
$$u_\e=V_{y_\e}+\phi_{y_\e}\in \mathcal H_{s}, \ \ \ \ y_\e=(r_\e,d_\e), \ \ \ \ \mbox{with} \ \ \|\phi_{y_\e}\|\leqslant C\vartheta(\mu).$$ Moreover, up to a subsequence, $r_\e \to r_0$ and $d_\e \to d_0$, where $r_0$ is a global minimum point of $\Phi$ and
\[d_0=\left(\frac{2B}{(N-2)A\,\Phi(r_0)}\right)^{\frac1{N-4}}.\]
\end{proposition}

\begin{proof}
By \eqref{hypPhi}, the set $\mathcal K$ of global minimum points of $\Phi$ is a nonempty compact subset of $I$. For fixed $r$, the function $d\mapsto\Gamma(r,d)$ has a unique minimum point at $d(r)$, and by \eqref{condminphi}, we can consider
\[
f(\Phi(r)) : = \Gamma(r,d(r))=-\frac{N-4}{2(N-2)}\left(\frac2{(N-2)\Phi(r)}\right)^{\frac2{N-4}},
\]
which is an increasing function of $\Phi(r)$. Hence $\Gamma$ attains its global minimum on $\mathcal M$ exactly on the compact set $\{(r,d(r)):\ r\in\mathcal K\}$. Choosing 
$$r_-<\min\mathcal K\leqslant\max\mathcal K<r_+$$
in $I$ and $0<d_-<d(r)<d_+$ for $r\in[r_-,r_+]$, the box 
$$\mathcal C:=[r_-,r_+]\times[d_-,d_+]$$satisfies
\[
\min_{\partial\mathcal C}\Gamma>\min_{\mathcal C}\Gamma.
\]
Fix $\eta>0$ such that $\{ (r, cd) : (r,d) \in \mathcal C\} \subset \mbox{int} \ P_{\eta}$. By Proposition~\ref{prop:expansion}, the function
\[
F_\e(r,d):=\frac{\mu^{-(N-2)}}{\ell Bc^2}\left(\widetilde J_\e(r,cd)-\ell\frac{S^{N/2}}{N}\right)
\]
converges to $\Gamma$ uniformly on $\mathcal C$. Hence, for $\e$ small, $\min_{\partial\mathcal C}F_\e>\min_{\mathcal C}F_\e$ and $F_\e$ has a minimum point $(r_\e,\bar d_\e)$ in the interior of $\mathcal C$. Then $y_\e:=(r_\e,c\bar d_\e)$ is a critical point of $\widetilde J_\e$ and, by Lemma~\ref{lem:reduced} and Proposition~\ref{prop:reduction}, we can conclude that $u_\e:=V_{y_\e}+\phi_{y_\e}$ is a solution of \eqref{pb} with the stated estimate. By uniform convergence, every limit point $(r_0,\bar d_0)$ of $(r_\e,\bar d_\e)$ is a minimum point of $\Gamma$ on $\mathcal C$, i.e.\ $r_0\in\mathcal K$ and $\bar d_0=d(r_0)$, which gives the value of $d_0=c\,d(r_0)$. Finally, $u_\e\not\equiv0$ and $u_\e\circ R=-u_\e$ in case (A), $u_\e\circ T=-u_\e$ in case (B), so that $u_\e$ changes sign and this completes the proof.
\end{proof}
\begin{remark}\label{rem:abstract}
By Lemma~\ref{lem:critPsi}, any critical point of $\Phi$ at which $\Phi>0$ gives a critical point of $\Gamma$ and a non-degenerate one gives a non-degenerate critical point. Since $\varpi_\e\to0$ in $C^1$, a standard degree argument (as in \cite[Section 3.1.3]{PRV3}) gives a solution concentrating at the corresponding configuration. In this paper we only use global minima, which are $C^0$-stable and do not require any non-degeneracy.
\end{remark}

\section{The case of the annulus}\label{anello}
\noindent Let us consider  the annulus $\Omega_{\rho}$ as in \eqref{omegarho} and $k\ge2$ be an even integer. \ \\ 
We are in case $(A)$ of Subsection \ref{subsec:sym}, so that the function $\Phi(r)$ in \eqref{Phigeneral} is
\begin{equation}\label{Psianello}
\Psi_k(r):=\tau_{\Omega_\rho}(\xi_1(r))-\sum_{j=1}^{k-1}(-1)^jG_{\Omega_\rho}(\xi_1(r),\xi_{j+1}(r)),\qquad r\in(\rho,1).
\end{equation}
By Proposition~\ref{propabstract}, Theorem~\ref{main}\,(a) follows once we show that $\Psi_k>0$ on $(\rho,1)$ and $\Psi_k(r)\to+\infty$ as $r\to\rho^+$ and as $r\to1^-$.
It is standard the behavior of the function $\Psi_{k}$ at the endpoints, instead it is much more involved the proof of $\Psi_k>0$. \ \\
So, in order to show that, we write \eqref{Psianello} by using the explicit expression of the Robin function given by (see \cite{GV}) 
\begin{equation*}
\tau_{\Omega_\rho}(x)=\frac{1}{\omega_{N-1}}\sum_{m=0}^{\infty}d_m Q_m(|x|),
\end{equation*}
where, for all $m \geq 0$,  
\begin{equation*}
d_{m} := \binom{m+N-2}{N-2} + \binom{m+N-3}{N-2}
\end{equation*}
and
\begin{equation}\label{Qm}
Q_m(|x|):=\frac{\rho^{2m+N-2} - 2\rho^{2m+N-2}|x|^{2m+N-2} + |x|^{4m+2N-4}}{(2m+N-2)|x|^{2N+2m-4}(1-\rho^{2m+N-2})}.
\end{equation}
We point out that,
$$\rho^{2m+N-2}-2\rho^{2m+N-2}r^{2m+N-2}+r^{2(2m+N-2)}=(r^{2m+N-2}-\rho^{2m+N-2})^2+\rho^{2m+N-2}(1-\rho^{2m+N-2})$$
and so, we have that $Q_m(r)>0$ for all $m\geqslant0$ and $r\in(\rho,1)$. \ \\
Moreover, the Green function of the annulus is given by 
\begin{equation*}
G_{\Omega_\rho}(x, y)=\frac{1}{\omega_{N-1}} \Bigg[ \frac{1}{(N-2)|x-y|^{N-2}} - \sum_{m=0}^{\infty} Q_m(x, y) Z_m\left(\frac{x}{|x|}, \frac{y}{|y|}\right) \Bigg],
\end{equation*}
where
\begin{equation*}
Q_m(x, y):=\frac{\rho^{2m+N-2} - \rho^{2m+N-2}(|x|^{2m+N-2}+|y|^{2m+N-2}) + |x|^{2m+N-2}|y|^{2m+N-2}}{(2m+N-2)(|x||y|)^{m+N-2}(1-\rho^{2m+N-2})},
\end{equation*}
where $\omega_{N-1} = \mbox{meas}(\mathbb{S}^{N-1})$and $Z_m(\zeta, \eta)$ represents the zonal harmonics of degree $m$ which have a particularly simple expression in terms of the Gegenbauer (or ultra-spherical) polynomials $C_{m}^{\lambda}$. The latter can be defined in terms of generating functions. If we write (see \cite{SW} p. 148)
\begin{equation}\label{generatric}
(1-2st+s^{2})^{-\lambda} = \sum_{m=0}^{\infty} C_m^{\lambda}(t)s^{m},
\end{equation}
where $0 \leq s<1$, $|t| \leq 1$ and $\lambda>0$, then the coefficient $C_m^{\lambda}$ is called Gegenbauer polynomial of degree $m$ associated with $\lambda$. \ \\
We have, for all $m \geq 0$,
\begin{equation}\label{binom}
d_{m} = A_{N,m} \cdot c_{N,m}\,,
\end{equation}
where 
\begin{equation}\label{ANm}
A_{N,m} : =\binom{N+m-3}{N-3}, \quad c_{N,m} = \frac{N+2m-2}{N-2}.
\end{equation}
An important fact is to understand the link between the zonal harmonics and the Gegenbauer polynomials. Indeed (see Theorem 2.1 in \cite{GV})
%\begin{theorem}[\cite{GV}, Theorem $2.1$]\label{zmpm}
if $N>2$ and $m=0, 1, 2, \ldots$ then we have that for all $x',y'$ such that $|x'|=|y'|=1$ it holds
\begin{equation}\label{zmcm}
Z_{m}(x',y')= c_{N,m} C_{m}^{\frac{N-2}{2}}(x' \cdot y').
\end{equation}
Inserting these expressions in \eqref{Psianello}, we have to study the sign of
\begin{equation}\label{expaut}
\Psi_k(r)= \frac{1}{\omega_{N-1}} \Bigg[ \sum_{m=0}^{+\infty} d_{m} Q_{m}(r) - \sum_{j=1}^{k-1} (-1)^{j} \Bigg( \frac{1}{(N-2)|\xi_1-\xi_{j+1}|^{N-2}} - \sum_{m=0}^{+\infty} Q_{m}(r) Z_{m} \left(\frac{\xi_1}{|\xi_1|}, \frac{\xi_{j+1}}{|\xi_{j+1}|} \right) \Bigg) \Bigg].
\end{equation}
%and so
%\begin{equation}\label{expaut}
%\Lambda_{0}(r) \omega_{N-1} = \Bigg[ \sum_{m=0}^{+\infty} d_{m} Q_{m}(\rho) - \sum_{j=1}^{k-1} (-1)^{j} \Bigg( \frac{1}{(N-2)|\xi_1-\xi_{j+1}|^{N-2}} - \sum_{m=0}^{+\infty} Q_{m}(\rho) Z_{m} \left(\frac{\xi_1}{|\xi_1|}, \frac{\xi_{j+1}}{|\xi_{j+1}|} \right) \Bigg) \Bigg]
%\end{equation}
Firstly we point out that
\begin{equation}\label{xi}
\begin{aligned}
|\xi_1(r) - \xi_{j+1}(r)|^{N-2} & = r^{N-2} \Bigg(2\sin \left( \frac{\pi j}{k} \right)\Bigg)^{N-2}\qquad \mbox{for all}\, j=1, \cdots, k-1.
%\left( |\xi_1(\rho) - \xi_{j+1}(\rho)|^{2} \right)^{\frac{N-2}{2}} = \left( 2r^2 - 2 \xi_1(\rho) \cdot \xi_{j+1}(\rho) \right)^{\frac{N-2}{2}} 
%\\& = \Bigg( 2\rho^2 - 2\rho^2 \cos \left( \frac{2\pi j}{k} \right) \Bigg)^{\frac{N-2}{2}} = \rho^{N-2} \Bigg(2\sin \left( \frac{\pi j}{k} \right)\Bigg)^{N-2}.
\end{aligned}
\end{equation}
 
Hence
\begin{equation*}
\sum_{j=1}^{k-1} \frac{(-1)^{j}}{|\xi_1(r)-\xi_{j+1}(r)|^{N-2}} = \frac{\mathfrak C_{k,N}}{r^{N-2}},\qquad \mbox{where}\,\,\, \mathfrak C_{k,N} : =\sum_{j=1}^{k-1} \frac{(-1)^{j}}{\left( 2 \sin \left( \frac{\pi j}{k} \right) \right)^{N-2}}.
\end{equation*}
%where 
%\begin{equation*}
%\mathfrak C_{k,N} : =\sum_{j=1}^{k-1} \frac{(-1)^{j}}{\left( 2 \sin \left( \frac{\pi j}{k} \right) \right)^{N-2}}
%\end{equation*} 
We can rewrite \eqref{expaut} as 
\begin{equation}\label{autov2}
\Psi_k(r) \omega_{N-1} = -\frac{\mathfrak C_{k,N}}{(N-2)r^{N-2}} + \sum_{m=0}^{+\infty} \Bigg[d_m+ \sum_{j=1}^{k-1}(-1)^{j} Z_{m} \left( \frac{\xi_1}{|\xi_1|}, \frac{\xi_{j+1}}{|\xi_{j+1}|}\right) \Bigg] Q_{m}(r).
\end{equation}
Now, using \eqref{zmcm} and \eqref{binom}, \eqref{autov2} becomes
\begin{equation}\label{aim10}
\Psi_k(r) \omega_{N-1} = -\frac{\mathfrak C_{k,N}}{(N-2)r^{N-2}} + \sum_{m=0}^{+\infty} c_{N,m} \alpha_{k,m,N} Q_{m}(r),
\end{equation}
where 
\begin{equation}\label{alpha}
\alpha_{k,m,N} : = A_{N,m} + S_{k,m,N} \ \ \ \mbox{with} \ \ S_{k,m,N} : = \sum_{j=1}^{k-1} (-1)^{j} C_{m}^{\lambda} \left( \cos \left( \frac{2\pi j}{k} \right) \right)
\end{equation}
and $\lambda = \frac{N-2}{2}$. \ \\
Now, since $Q_{m}(r)>0$ for all $m\geq 0$ and for any $\rho<r<1$, we have to show that \begin{equation}\label{aim12}
\mathfrak C_{k,N} < 0  \ \ \ \mbox{and} \ \ \alpha_{k,m,N} \geq 0 \,\,\mbox{for all}\,\, m\geq 0.
\end{equation}
\noindent {\it Aim 1:} $\mathfrak C_{k,N} < 0$.\\
\noindent
First, since $ \sin \left( \frac{\pi j}{k} \right) = \sin \left( \frac{\pi (k-j)}{k} \right)$ and $(-1)^{k-j} = (-1)^{j}$, we can write 
%It is easily to observe that, since $ \sin \left( \frac{\pi j}{k} \right) = \sin \left( \frac{\pi (k-j)}{k} \right)$ and $(-1)^{k-j} = (-1)^{j}$ (since $k$ is even), 
\begin{equation*}
\mathfrak C_{k,N} = \frac{1}{2^{N-2}} \left( 2\sum_{j=1}^{\frac{k}{2}-1} \frac{(-1)^{j}}{\sin^{N-2} \left( \frac{\pi j}{k} \right)} + (-1)^{\frac{k}{2}}\right).
\end{equation*}

If $\frac{k}{2}$ is odd, we have that
\begin{equation*}
\begin{aligned}
\mathfrak C_{k,N} & = \frac{1}{2^{N-2}} \left( 2\sum_{j=1}^{\frac{k}{2}-1} \frac{(-1)^{j}}{\sin^{N-2} \left( \frac{\pi j}{k} \right)} -1\right) 
\\& = \frac{1}{2^{N-2}} \left(-1+2\sum_{h=1}^{\frac{k}{4}-\frac{1}{2}} -\frac{1}{\sin^{N-2} \left(\frac{\pi}{k}(2h-1)\right)} + 2\sum_{h=1}^{\frac{k}{4}-\frac{1}{2}} \frac{1}{\sin^{N-2} \left(\frac{2\pi h}{k}\right)} \right)
\\& = \frac{1}{2^{N-2}} \Bigg[ -1 + 2 \sum_{h=1}^{\frac{k}{4}-\frac{1}{2}} \left( \frac{1}{\sin^{N-2} \left(\frac{2\pi h}{k}\right)} - \frac{1}{\sin^{N-2} \left(\frac{\pi}{k}(2h-1)\right)} \right) \Bigg].
\end{aligned}
\end{equation*}
Now, we can observe that, since $h \leq \frac{k}{4}$,  
$$0<\frac{\pi}{k}(2h-1)<\frac{2 \pi h}{k} < \frac{\pi}{2}$$
and so, from the fact that $\sin x$ is increasing in $(0, \frac{\pi}{2})$ we observe that 
\begin{equation*}
\frac{1}{\sin^{N-2} \left(\frac{\pi}{k}(2h-1)\right)} > \frac{1}{\sin^{N-2} \left(\frac{2\pi h}{k}\right)}.
\end{equation*}
Hence, 
\begin{equation*}
\mathfrak C_{k,N} =\frac{1}{2^{N-2}}\left[ -1 +  2 \sum_{h=1}^{\frac{k}{4}-\frac{1}{2}} \left( \frac{1}{\sin^{N-2}\left(\frac{2\pi h}{k}\right)} - \frac{1}{\sin^{N-2} \left(\frac{\pi}{k}(2h-1)\right)} \right) \right]< 0 
\end{equation*}
and so the claim follows. \ \\
Now, we want to consider the case in which $\frac{k}{2}$ is even. Then we have that
\begin{equation*}
\mathfrak C_{k,N} = \frac{1}{2^{N-2}} (1 + 2 \mathcal A), \ \ \ \mbox{with} \ \ \mathcal A : = \sum_{j=1}^{\frac{k}{2}-1} \frac{(-1)^{j}}{\sin^{N-2} \left(\frac{\pi j}{k} \right)}
\end{equation*}
and, so, in order to get $\mathfrak C_{k,N} < 0$, we need to prove that $\mathcal A < -\frac12$. \ \\
Now, we can write 
\begin{equation*}
\mathcal A = \sum_{j=1}^{\frac{k}{2}-1} (-1)^{j} a_{j}, \ \ \ \mbox{where} \ \ a_{j} : = \frac{1}{\sin^{N-2}\left( \frac{\pi j}{k} \right)}.
\end{equation*}
Now, since $\frac{k}{2}-1$ is odd, we have that
$$\mathcal A = - (a_1 -a_2) - (a_3 - a_4 ) - \cdots - (a_{\frac{k}{2}-3} - a_{\frac{k}{2}-2}) - a_{\frac{k}{2}-1} < - a_{\frac{k}{2}-1}$$
since $a_{j}>a_{j+1}$. \ \\
Moreover, we can observe that
\begin{equation*}
a_{\frac{k}{2}-1} = \frac{1}{\sin^{N-2} \left(\frac{\pi}{k}\left(\frac{k}{2}-1\right)\right)} = \frac{1}{\cos^{N-2} \left(\frac{\pi}{k}\right)} >1.
\end{equation*}\ \\
Therefore, we obtain that $\mathcal A < -1$.
Hence the thesis follows. \\\\ \\
\noindent
{\it Aim 2:} $\alpha_{k,m,N}\geq 0$ for all $N \geq 3$, $k \geq 2$ even and $m\geq 0$. \ \\ \\
Now, letting $\beta_{j} : = \frac{2\pi j}{k}$, we recall the following expression of Gegenbauer polynomials (see \cite[8.934.2]{Gra})
\[
C^\lambda_m(\cos\beta_{j})=\sum_{l=0}^m\gamma_{m,l}\cos\bigl((m-2l)\beta_{j}\bigr),\qquad \gamma_{m,l}:=\frac{\Gamma(\lambda+l)\,\Gamma(\lambda+m-l)}{\Gamma(\lambda)^2\,l!\,(m-l)!}>0 .
\]
Taking $t=1$ in \eqref{generatric} and using the binomial series, for $|s|<1$
\[
\sum_{m=0}^\infty C^\lambda_m(1)\,s^m=(1-s)^{-2\lambda}=\sum_{m=0}^\infty\frac{\Gamma(2\lambda+m)}{\Gamma(2\lambda)\,m!}\,s^m ,
\]
so that, by the uniqueness of the coefficients of a power series and since $2\lambda=N-2$,
\[
C^\lambda_m(1)=\frac{(N+m-3)!}{(N-3)!\,m!}=A_{N,m}.
\]
Moreover, taking $\beta_j=0$ in the expression of $C^\lambda_m(\cos\beta_j)$ above, $\displaystyle{C^\lambda_m(1)=\sum_{l=0}^m\gamma_{m,l}}$. \ \\
Therefore 
$$A_{N,m} =C^\lambda_m(1)= \sum_{l=0}^m\gamma_{m,l};$$
so, since $(-1)^j\cos((m-2l)\beta_{j})=\cos((m-2l)\beta_{j}+\pi j)$ and $\pi j=\frac k2\cdot\frac{2\pi j}k$,
\[
\alpha_{k,m,N}=\sum_{l=0}^m\gamma_{m,l}\Biggl[1+\sum_{j=1}^{k-1}\cos\left(\left(m-2l+\frac k2\right)\frac{2\pi j}k\right)\Biggr].
\]
In addition, for every integer $n$ we have (see also \cite[(4.29)]{PRV2})
 \[
\sum_{j=0}^{k-1} e^{2\pi i n j/k}
=
\begin{cases}
k & \text{if } k \mid n,\\
0 & \text{otherwise};
\end{cases}
\]hence
\begin{equation}\label{portantemk}
1+\sum_{j=1}^{k-1}\cos\Bigl(\frac{2\pi nj}k\Bigr)=\begin{cases}k & \text{if $k$ divides $n$},\\ 0 & \text{otherwise}.\end{cases}
\end{equation}
Applying \eqref{portantemk} with the integer $n=m-2l+\frac k2$ we can conclude that
\[
\alpha_{k,m,N}=k\sum_{\substack{0\leqslant l\leqslant m\\ k\,\mid\, m-2l+k/2}}\gamma_{m,l}\ \geqslant0.
\]
By \eqref{aim10} and \eqref{aim12}, 
$$\omega_{N-1}\Psi_k(r)\geqslant-\frac{\mathfrak C_{k,N}}{(N-2)r^{N-2}}>0$$ \ \\
for every $r\in(\rho,1)$, which concludes the proof of Theorem \ref{main}\,(a).

\section{The case of the ball}\label{palla}
\noindent Let us consider, for every $N \geq 5$, the unit ball
\begin{equation*}
\mathbb{B}_1(0) : = \{ x \in \mathbb{R}^{N} : |x|<1 \}.
\end{equation*} 
The Robin function and the Green function of $\mathbb{B}_{1}(0)$ are explicitly given by
$$\tau(x):=\frac{1}{(1-|x|^2)^{N-2}},\quad 
G(x, y):=\left(\frac{1}{|x-y|^{N-2}}-\frac{1}{(1-2 \langle x,y \rangle + |x|^{2}|y|^{2})^{\frac{N-2}{2}}} \right)\\$$ 
up to the multiplicative constant $\omega_{N}^{-1}$, which we drop in order to simplify the computations (it can be absorbed in the constants $A,B$ of Proposition \ref{prop:expansion}). \ \\
\subsection{The symmetric setting}
We are in case $(B)$ of Subsection \ref{subsec:sym}, so that the function \eqref{Phigeneral} is 
\begin{equation}\label{eq:Phi-def}
\Phi_{k,N}(r):=\tau(\xi_1)-\sum_{j=2}^{k}G(\xi_1,\xi_j)+\sum_{j=1}^{k}G(\xi_1,\eta_j),
\qquad r\in(0,1),
\end{equation}
where $\xi_j=\xi_j(r)$ and $\eta_j=\eta_j(r)$.  \ \\ \\
By Proposition \ref{propabstract}, Theorem \ref{main} $(b)$ follows once we prove that $\Phi_{k,N}$ satisfies \eqref{hypPhi} on $I=(0,1)$ for $2 \leq k \leq k^{*}(N)$. \ \\ \\
Now, we let
\begin{equation}\label{lambdak}
\Lambda_{k} : = \sum_{j=1}^{k-1}\frac{1}{\sin^{N-2}\Big(\frac{\pi j}{k}\Big)}
\end{equation} 
and we introduce the following constant
\begin{equation}\label{eq:ckN}
\mathfrak c_{k,N}:=k\,2^{-\frac{(N-2)}{2}}-2^{-(N-2)}\Lambda_{k},
\end{equation}
or, equivalently, $D_k(N)=2^{N-2} \mathfrak c_{k,N}=k2^{\frac{N-2}{2}}-\Lambda_{k}$.

\subsubsection{Structure of $\Phi_{k,N}$}\label{sec:Phi}

In order to evaluate the Green and the Robin function we recall some useful facts. \ \\
In particular, for $j=1,\dots,k$,
\[|\xi_1-\xi_j|^2=2r^2\big(1-\cos\theta_j\big)=4r^2\sin^2\frac{\pi(j-1)}{k},\qquad\xi_1\cdot\xi_j=r^2\cos\theta_j,\]
while, since the two planes are orthogonal,
\[|\xi_1-\eta_j|^2=2r^2,\qquad \xi_1\cdot\eta_j=0 .\]
Fix $r \in(0,1)$ and define
\begin{equation}\label{eq:h}
h(t)=h_r(t):=\big(1+r^4-2r^2t\big)^{-\frac{(N-2)}{2}},\qquad t\in[-1,1].
\end{equation}
Since $1+r^4-2r^2t\geq(1-r^2)^2>0$ on $[-1,1]$, $h$ is smooth and positive there, and
\begin{equation}\label{notimp}
H(\xi_1,\xi_j)=h(\cos\theta_j),\qquad \tau(\xi_1)=h(1)=(1-r^2)^{-(N-2)},\qquad H(\xi_1,\eta_j)=h(0)=(1+r^4)^{-\frac{(N-2)}{2}}.
\end{equation}
\begin{lemma}[Decomposition]\label{lem:decomp}
For every $k\geq2$ and $r \in(0,1)$,
\begin{equation}\label{eq:decomp}
\Phi_{k,N}(r)=R_{k,N}(r)+\frac{\mathfrak c_{k,N}}{r^{N-2}},
\qquad
R_{k,N}(r):=\sum_{j=1}^{k}h(\cos\theta_j)-k\,h(0),
\end{equation}
with $\mathfrak c_{k,N}$ given by~\eqref{eq:ckN}.
\end{lemma}
\begin{proof}
In the following we will treat separately the regular and the singular part of the Green function. In particular, 
\begin{equation}
\begin{aligned}
\Phi_{k,N}(r) & = \tau(\xi_{1}) - \sum_{j=2}^{k} \Big[ \frac{1}{|\xi_1-\xi_{j}|^{N-2}} - H(\xi_1, \xi_j) \Big]+ \sum_{j=1}^{k} \Big[ \frac{1}{|\xi_1-\eta_j|^{N-2}} - H(\xi_1, \eta_j) \Big] 
\\& =\underbrace{\sum_{j=1}^{k} H(\xi_1, \xi_j) - \sum_{j=1}^{k} H(\xi_1, \eta_j)}_{(I_1)} + \underbrace{\sum_{j=1}^{k} \frac{1}{|\xi_1-\eta_j|^{N-2}}-\sum_{j=2}^{k} \frac{1}{|\xi_1-\xi_j|^{N-2}}}_{(I2)}.
\end{aligned}
\end{equation}
Therefore, from \eqref{notimp}, the regular parts give
\[(I_1)=\tau(\xi_1)+\sum_{j=2}^kH(\xi_1,\xi_j)-\sum_{j=1}^kH(\xi_1,\eta_j)=h(\cos\theta_1)+\sum_{j=2}^kh(\cos\theta_j)-k\,h(0)=R_{k,N}(r).\]
Instead,  the singular parts give
\[
\begin{aligned}
(I_{2})
&= - \sum_{j=2}^{k}
\frac{1}{\left(2r\sin\frac{\pi(j-1)}{k}\right)^{N-2}}
+ k(2r^2)^{-\frac{N-2}{2}}
\\
&= r^{-(N-2)}
\left[
k\,2^{-\frac{N-2}{2}}
-2^{-(N-2)}
\sum_{j=1}^{k-1}
\frac{1}{\sin^{N-2}\left(\frac{\pi j}{k}\right)}
\right]
= \frac{\mathfrak c_{k,N}}{r^{N-2}}.
\qedhere
\end{aligned}
\]
\end{proof}
We now study the sign of $\Phi_{k,N}$, analyzing separately the sign of the regular part and of the singular part. \ \\ \\
%\textcolor{red}{Giuseppe: Quello che si vuole fare ora è legato allo studio del segno del funzionale utilizzando la decomposizione precedente, si parte dall'analisi della positività della parte regolare e di quella singolare successivamente. In più nel Lemma che segue si ottiene una buona stima dal basso della parte regolare necessaria nel momento in cui, per assicurarci l'esistenza del minimo, andremo a studiare i limiti agli estremi del dominio della $\Phi_{k,N}$. 
%Tuttavia non so se sia il modo migliore di scriverlo o si possa scrivere meglio, tu che dici?} \ \\
\begin{lemma}[Positivity of the regular part]\label{lem:R}
For every $k\geq2$ and every $r\in(0,1)$,
\begin{equation}\label{eq:R-bounds}
R_{k,N}(r)>\max\Big\{0,\ (1-r^2)^{-(N-2)}-k\Big\}.
\end{equation}
In particular $R_{k,N}(r)\to+\infty$ as $r\to1^-$.
\end{lemma}
\begin{proof} 
We have
\begin{equation*}
h'(t)=(N-2) r^2(1+ r^4-2r^2t)^{-\frac{N}{2}}>0
\end{equation*} 
and
\[
h''(t)=N(N-2)r^4\big(1+r^4-2r^2t\big)^{-\frac{N}{2}-1}>0\qquad\text{on }[-1,1],
\]
so $h$ is strictly convex. \ \\
Furthermore, it is known that 
%\textcolor{red}{ho controllato, è vera e mi sembra sia una identità standard}
\[\frac1k\sum_{j=1}^{k}\cos\frac{2\pi(j-1)}{k}=\frac1k\,\mathrm{Re}\sum_{j=0}^{k-1}e^{2\pi ij/k}=0 .\]
By the strict Jensen inequality,
\[
\frac1k\sum_{j=1}^k h(\cos \theta_j)>h\left(\frac1k\sum_{j=1}^k \cos \theta_j\right)=h(0)
\]
and so, we can conclude that $R_{k,N}(r)>0$. \ \\ \\
Finally, since $h>0$ and $h(0)=(1+r^4)^{-\frac{(N-2)}{2}} <1$,
\begin{equation*}
\begin{aligned}
R_{k,N}(r) & = \sum_{j=1}^{k} h(\cos \theta_{j}) -kh(0) = h(\cos \theta_1)+\sum_{j=2}^{k} h(\cos \theta_{j}) - kh(0) 
\\& = h(1) + \sum_{j=2}^{k} \underbrace{h(\cos \theta_{j})}_{>0} - kh(0) > h(1)-kh(0) \geq (1-r^{2})^{-(N-2)} -k. 
\end{aligned}
\end{equation*}
We also note that $h_{r}(t) \to1$ uniformly on $[-1,1]$ as $r\to 0^{+}$, so that $R_{k,N}(r) \to 0$ as $r \to 0^{+}$.
\end{proof}
\ \\
%\textcolor{blue}{ Giuseppe : Al solito, se la funzione responsabile dell'esistenza del punto critico che bramiamo è interamente positiva e va a $+\infty$ da ambo gli estremi della configurazione, dovrà necessariamente avere un minimo. \ \\
%Per quanto riguarda la $(iii)$ mi sembra serva per dire che $\Phi_{k,N}$ ha un minimo globale positivo e che l'insieme dei punti di minimo sta in un compatto di $(0,1)$, lontano sia da $0$ sia da $1$.}
\begin{proposition}\label{prop:Phi-min}
Let $k\geq2$ and assume $\mathfrak c_{k,N}>0$. Then:
\begin{enumerate}[label=(\roman*)]
\item $\Phi_{k,N}(r)>\mathfrak c_{k,N}r^{-(N-2)}>0$ for every $r\in(0,1)$;
\item $\Phi_{k,N}(r)\to+\infty$ as $r\to0^+$ and as $r\to1^-$.
\end{enumerate}
In particular, $\Phi_{k,N}$ satisfies \eqref{hypPhi} on $I=(0,1)$.
\end{proposition}
\begin{proof}
We can observe that $(i)$ follows from Lemmas \ref{lem:decomp} and \ref{lem:R}. \ \\
Instead, we get that, as $r \to 0^{+}$,  
$$\Phi_{k,N}(r)\geq \frac{\mathfrak c_{k,N}}{r^{N-2}} \to +\infty$$ 
and, as $r \to1^-$, $\Phi_{k,N}(r)\geq R_{k,N}(r)\to+\infty$. Therefore $(ii)$ easily follows. 
\end{proof}

%\begin{remark}[What happens if $\mathfrak c_{k,N}<0$]\label{rem:negative}
%Since $R_{k,N}$ is bounded near $0$, if $\mathfrak c_{k,N}<0$ then $\Phi_{k,N}(r)\to-\infty$ as $r\to0^+$, and $\Phi_{k,N}<0$ on some interval $(0,r_1)$. By Lemma \ref{lem:critPsi}, the reduced function $\Gamma$ has no critical points with $r$ in this interval and the global-minimum argument of Proposition \ref{propabstract} breaks down. This does not prove that no solution of this shape exists, but that this method gives nothing. Numerically, for the values $(k,N)$ with $\mathfrak c_{k,N}< 0$ that we checked, $\Phi_{k,N}$ is strictly increasing on $(0,1)$, exactly as for $k\geq 7$ in the $4D$-case (\cite{PRV2}).
%\end{remark}

\subsubsection{The sign of $\mathfrak c_{k,N}$}\label{sec:sign}

\begin{proposition}\label{prop:sign}
Let $N \geq 5$ and $D_k(N)=2^{N-2}\mathfrak c_{k,N}$ as in~\eqref{eq:ckN}. Then
\begin{enumerate}[label=(\alph*)]
\item $D_k(N)>0$ for $k=2,3,4$ and every $N\geq5$;
\item $D_5(5)>0$ and $D_5(6)=\tfrac45>0$, while $D_5(N)<0$ for every $N\geq7$;
\item $D_k(N)<0$ for every $k\geq6$ and every $N\geq5$.
\end{enumerate}
Consequently, for $N\geq5$ and $k \geq 2$, $\mathfrak c_{k,N}>0$ if and only if $k\leq k^*(N)$, with $k^*(N)$ as in Theorem \ref{main}.
\end{proposition}

\begin{proof}
$(a)$ A direct computation gives 
\begin{itemize}
\item $k=2$: $D_2(N)=2^{\frac{N}{2}}-1>0$;
\item $k=3$: $D_{3}(N) = \left( \frac{2}{3} \right)^{\frac{N-2}{2}} (3^{\frac{N}{2}}-2^{\frac{N}{2}}) >0$;
\item $k=4$: $D_{4}(N) = 2^{\frac{N}{2}}-1>0$.
\end{itemize} \ \\
$(b)$ We use the classical values
\[
\sin^2\left(\frac\pi5\right)=\frac{5-\sqrt5}{8} \ \ \ \mbox{and} \ \ \ \ \sin^2\left(\frac{2\pi}5 \right)=\frac{5+\sqrt5}{8},
\]
and we recall that $\sin\left(\frac{3\pi}{5}\right)=\sin\left(\frac{2\pi}{5}\right) \ \ \ \mbox{and} \ \ \ \ \sin\left(\frac{4\pi}{5}\right)=\sin\left(\frac{\pi}{5}\right).$ \ \\
Set 
\[
\mathfrak a : = \left(2\sin^2\left(\frac{\pi}{5}\right)\right)^{-1}=\frac{4}{5-\sqrt{5}} \ \ \ \ \mbox{and} \ \ \ \ 
\mathfrak b : = \left(2\sin^2\left(\frac{2\pi}{5}\right)\right)^{-1}=\frac{4}{5+\sqrt{5}};
\]
hence 
$$\sin^{-(N-2)}\left(\frac{\pi j}{5}\right)=2^{\frac{N-2}{2}}\mathfrak a^{\frac{N-2}{2}} \ \ \mbox{for} \ \ \ j=1,4$$
and 
$$\sin^{-(N-2)}\left(\frac{\pi j}{5}\right)=2^{\frac{N-2}{2}}\mathfrak b^{\frac{N-2}{2}} \ \ \ \mbox{for} \ \ \ j=2,3.$$ 
Therefore
$$D_{5}(N) = 2^{\frac{N-2}{2}} \left[5-2(\mathfrak a^{\frac{N-2}{2}} + \mathfrak b^{\frac{N-2}{2}} )\right].$$
Note that 
$$\mathfrak a^{2} + \mathfrak b^{2} = 2\left(1+\frac{1}{5}\right) = \frac{12}{5}.$$
Moreover, for $p>0$, let
$$M_{p} := \left(\frac{\mathfrak a^{p}+\mathfrak b^{p}}{2} \right)^{\frac{1}{p}}$$
be the power mean of order $p$ of $\mathfrak a, \mathfrak b$; it is well known that $p \to M_{p}$ is nondecreasing and $\mathfrak a^{p}+\mathfrak b^{p} = 2M_{p}^{p}$ and in particular $M_{2} = \left(\frac{6}{5}\right)^{\frac{1}{2}}$. 
Therefore, it is easy to check that
$$D_{5}(6) = 4\left[5-2(\mathfrak a^{2}+\mathfrak b^{2})\right] = \frac{4}{5}.$$
Instead, since $M_{\frac{3}{2}} \leq M_{2}$, we get that 
\begin{equation*}
\mathfrak a ^{\frac{3}{2}} + \mathfrak b^{\frac{3}{2}} = 2M_{\frac{3}{2}}^{\frac{3}{2}} \leq 2M_{2}^{\frac{3}{2}} = 2 \left(\frac{6}{5}\right)^{\frac{3}{4}}.
\end{equation*}
Thus 
$$D_{5}(5) \geq 2^{\frac{3}{2}} \left[5-4\left(\frac{6}{5}\right)^{\frac{3}{4}}\right] >0.$$
Now, if we consider the case in which $N \geq 7$, using the fact that $\frac{N-2}{2} \geq \frac{5}{2} >2$ and so $M_{\frac{N-2}{2}} \geq M_{2} >1$, we have that 
$$\mathfrak a^{\frac{N-2}{2}} +\mathfrak b^{\frac{N-2}{2}}  =2M_{\frac{N-2}{2}}^{\frac{N-2}{2}} \geq 2M_{2}^{\frac{N-2}{2}} = 2 \left(\frac{6}{5} \right)^{\frac{N-2}{4}} \geq 2 \left(\frac{6}{5}\right)^{\frac{5}{4}}$$
and therefore
$$D_{5}(N) \leq 2^{\frac{N-2}{2}} \left[5-4\left(\frac{6}{5}\right)^{\frac{5}{4}}\right]<0.$$ \ \\ \\
$(c)$ We use the classical identity (see \cite{Fi})
\begin{equation}\label{eq:sin2}
\sum_{j=1}^{k-1}\frac{1}{\sin^2(\pi j/k)}=\frac{k^2-1}{3}.
\end{equation}
The numbers $y_j:=\sin^{-2}\left(\frac{\pi j}{k}\right)$, $j=1,\dots,k-1$, are positive and $y\mapsto y^{\frac{N-2}{2}}$ is convex on $(0,\infty)$ since $\frac{N-2}{2}\geq1$. By Jensen's inequality and~\eqref{eq:sin2},
\[\sum_{j=1}^{k-1}\frac{1}{\sin^{N-2}\left(\frac{\pi j}{k}\right)}=\sum_{j=1}^{k-1}y_j^{\,\frac{N-2}{2}}\geq(k-1)\left(\frac1{k-1}\sum_{j=1}^{k-1}y_j\right)^{\frac{N-2}{2}}=(k-1)\left(\frac{k+1}{3}\right)^{\frac{N-2}{2}}.\]
Therefore
\begin{equation}\label{DkN}
D_k(N)\leq k\,2^{\frac{N-2}{2}}-(k-1)\left(\frac{k+1}{3}\right)^{\frac{N-2}{2}}=k\,2^{\frac{N-2}{2}}\left[1-\frac{k-1}{k}\left(\frac{k+1}{6}\right)^{\frac{N-2}{2}}\right].
\end{equation}
Let now $k\geq6$. Then $\frac{k+1}{6}>1$, so 
$$\left(\frac{k+1}{6}\right)^{\frac{N-2}{2}}\geq\left(\frac{k+1}{6}\right)^{\frac{3}{2}}$$
for $N\geq5$; moreover $\frac{k-1}{k}$ and $\frac{k+1}{6}$ are increasing in $k$. \ \\
Hence
\begin{equation}\label{stimak}
\frac{k-1}{k}\Big(\frac{k+1}{6}\Big)^{\frac{N-2}{2}}\geq\frac56\Big(\frac76\Big)^{\frac{3}{2}}>1.
\end{equation}
Thus, using \eqref{DkN} and \eqref{stimak} we get that $D_k(N)<0$ for every $k\geq6$ and every $N\geq5$ and this concludes the proof.
\end{proof}

\begin{remark}\label{rmk4D}
For $N=4$ the identity~\eqref{eq:sin2} gives 
$$D_k(4)=2k-\frac{k^2-1}{3},$$
i.e. $\mathfrak c_{k,4}=\frac{-k^2+6k+1}{12}$, which is positive exactly for $k\leq6$: this is the coefficient appearing in the $4D$-case (see \cite{PRV2}).
\end{remark}
\begin{proof}[Proof of Theorem~\ref{main}\,(b)]
If $2\leq k\leq k^*(N)$, then $\mathfrak c_{k,N}>0$ by Proposition~\ref{prop:sign}; hence $\Phi_{k,N}$ satisfies \eqref{hypPhi} by Proposition~\ref{prop:Phi-min}, and Proposition~\ref{propabstract} applies in case (B). The convergence of the whole family follows from Corollary~\ref{cor:nondeg}.
\end{proof}

%\begin{remark}[Geometric meaning]
%In the energy, bubbles of the same sign \emph{attract} each other (the term $-G$ lowers $\Phi$ as they approach) and bubbles of opposite sign \emph{repel} (the term $+G$ raises $\Phi$). As $\rho\to0^+$ the attraction among the peaks of the same polygon, at distances $2\rho\sin(\pi j/k)$, tends to make the configuration collapse at the origin, and only the repulsion between the two polygons, at distance $\rho\sqrt2$, can prevent it: $c_{k,N}>0$ means exactly that the repulsion wins. As $N$ grows the interaction $|x-y|^{2-N}$ concentrates at short distances and only nearest neighbours matter: the repulsion wins if and only if $\rho\sqrt2<2\rho\sin(\pi/k)$, i.e. $k<4$, and $k=4$ is borderline (the side of the square equals the distance between the planes). This explains why $k^*(N)=4$ for large $N$.
%\end{remark}

%FORSE

%\subsubsection{Proof of Theorem \ref{thm:main} $(b)$}\label{sec:proof}
%Let $N \geq 5$ and $2\leq k \leq k_{*}(N). By Proposition \ref{prop:sign} $\mathfrak c_{k,N}>0$ and by Proposition \ref{prop: Phi-min}, $\Phi_{k,N}$ satisfies \eqref{HypPhi} on $I = (0,1)$. \ \\
%Therefore, Proposition \ref{propabstract}, applied in case $(B)$, gives a solution $u_{\e} \in \mathcal H_{s}$ of the form stated in Theorem 1.1 (b), with rε →r0, a
%1
%global minimum point of Φk,N, and dε →d0 =
%2B
%(N−2)AΦk,N(r0)
%N−4 . The symmetry properties
%are those of H; in particular uε(Tx) =−uε(x), so that uε changes sign.

\subsection*{The regular part as a power series}

Throughout this subsection $\lambda=\frac{N-2}{2}$, so that, with the notation of \eqref{eq:h},
\[
h_r(t)=\bigl(1-2r^2t+r^4\bigr)^{-\lambda}, \ \ \ \mbox{with} \ \ t\in[-1,1] \ \ \mbox{and} \ \ r\in(0,1).
\]
We recall from Section \ref{anello} the following expression of Gegenbauer polynomials, namely
\begin{equation}\label{gegfourier}
C^\lambda_m(\cos\varphi)=\sum_{l=0}^m\gamma_{m,l}\cos\bigl((m-2l)\varphi\bigr),\qquad
\gamma_{m,l}=\frac{\Gamma(\lambda+l)\,\Gamma(\lambda+m-l)}{\Gamma(\lambda)^2\,l!\,(m-l)!}>0
\end{equation}
and we set, for $n\geqslant0$,
\begin{equation}\label{bn}
b_n:=\frac{\Gamma(\lambda+n)}{\Gamma(\lambda)\,n!},\qquad\text{so that}\qquad \gamma_{2n,n}=b_n^2 .
\end{equation}

\begin{lemma}\label{lem:series}
Let $N\geqslant4$ and $k\geqslant2$. Then, for every $r\in(0,1)$,
\begin{equation}\label{Rseries}
R_{k,N}(r)=\sum_{m=1}^{\infty}a_m\,r^{2m},\qquad
a_m:=k\sum_{\substack{0\leqslant l\leqslant m\\ k\,\mid\,m-2l}}\gamma_{m,l}\;-\;k\,C^\lambda_m(0),
\end{equation}
where the series converges absolutely, $a_m\geqslant0$ for every $m\geqslant1$ and $a_2>0$.
\end{lemma}

\begin{proof}
By the generating function of the Gegenbauer polynomials, defined in \eqref{generatric}, with $r^{2}$ in place to $s$, we can write
\[
h_r(t)=\sum_{m=0}^\infty C^\lambda_m(t)\,r^{2m},\qquad t \in[-1,1].
\]
Moreover, by \eqref{gegfourier}, we get that
$$|C^\lambda_m(t)|\leqslant\sum_{l=0}^{m}\gamma_{m,l}=C^\lambda_m(1)$$
for $|t|\leqslant1$ and 
$$\sum_{m=0}^{+\infty} C^\lambda_m(1)\,r^{2m}=(1-r^2)^{-2\lambda}<\infty;$$
hence the series converges absolutely and uniformly in $t\in[-1,1]$. \ \\
Therefore
\[
R_{k,N}(r)=\sum_{j=1}^kh_r(\cos\theta_j)-k\,h_r(0)=\sum_{m=0}^\infty\Bigl(\sum_{j=1}^kC^\lambda_m(\cos\theta_j)-k\,C^\lambda_m(0)\Bigr)r^{2m}.
\]
By \eqref{gegfourier} and \eqref{portantemk} (recalling the definition of $\theta_j$),
\[
\sum_{j=1}^kC^\lambda_m(\cos\theta_j)=\sum_{l=0}^m\gamma_{m,l}\sum_{j=0}^{k-1}\cos\Bigl(\frac{2\pi(m-2l)j}{k}\Bigr)=k\sum_{\substack{0\leqslant l\leqslant m\\ k\,\mid\,m-2l}}\gamma_{m,l},
\]
which gives \eqref{Rseries}; the term $m=0$ vanishes, since $C^\lambda_0(t)\equiv1$ (see $8.930,1.10$ of \cite{Gra}). \ \\
Now, we have to compute $C^\lambda_m(0)$. Taking $t=0$ in the generating function of the Gegenbauer polynomials (see \cite[p.~148]{SW}) we get, for $|s|<1$,
\begin{equation}\label{gen0}
\sum_{m=0}^\infty C^\lambda_m(0)\,s^m=(1+s^2)^{-\lambda}.
\end{equation}
On the other hand, using standard results about the binomial series, we can write, for $|s|<1$,
\[
(1+s^2)^{-\lambda}=\sum_{m=0}^\infty c_m\,s^m,\qquad c_{2n}:=(-1)^n\,b_n,\quad c_{2n+1}:=0 .
\]
Both this series and the one in \eqref{gen0} are power series converging on $(-1,1)$ to the same function $f(s)=(1+s^2)^{-\lambda}$. Since the coefficients of a power series are determined by its sum, namely they are equal to $\frac{f^{(m)}(0)}{m!}$ and so we conclude that $C^\lambda_m(0)=c_m$ for every $m$, i.e.
\begin{equation}\label{C0}
C^\lambda_{2n+1}(0)=0,\qquad C^\lambda_{2n}(0)=(-1)^n\,b_n,\qquad n\geqslant0.
\end{equation}
Now we have to check that $a_m\geqslant0$. Since all the $\gamma_{m,l}$ are positive, the first term in the definition of $a_m$ is nonnegative.
\begin{itemize}[leftmargin=*]
\item If $m$ is odd, or $m=2n$ with $n$ odd, then $C^\lambda_m(0)\leqslant0$ and $a_m\geqslant0$.
\item If $m=2n$ with $n$ even, the index $l=n$ appears in the sum, because $m-2n=0$ is divisible by $k$. Hence
\[
\frac{a_{2n}}k\geqslant\gamma_{2n,n}-C^\lambda_{2n}(0)=b_n^2-b_n=b_n(b_n-1)\geqslant0,
\]
since, by the functional equations
$$\Gamma(x+1)=x\,\Gamma(x) \ \ \ \mbox{and} \ \ \ \Gamma(\lambda+n)=\Gamma(\lambda)\prod_{i=0}^{n-1}(\lambda+i),$$ 
easily follows
\[
b_n=\prod_{i=0}^{n-1}\frac{\lambda+i}{1+i}\geqslant1\qquad\text{as }\lambda\geqslant1 .
\]
\end{itemize}
Finally, with a straightforward computation we obtain that
$$a_2\geqslant k\bigl(\gamma_{2,1}-C^\lambda_2(0)\bigr)=k(\lambda^2+\lambda)>0.$$
\end{proof}
\begin{remark}
Lemma~\ref{lem:series} gives a second proof of the positivity $R_{k,N}>0$ in Lemma~\ref{lem:R}, and shows moreover that $R_{k,N}$ is strictly increasing on $(0,1)$.
\end{remark}

\subsection*{The shape of $\Phi_{k,N}$}

\begin{proposition}\label{prop:shape}
Let $N\geqslant5$ and $k\geqslant2$, and set
\[g(r):=r^{N-1}\,\Phi_{k,N}'(r),\qquad r\in(0,1).\]
Then $g$ is strictly increasing on $(0,1)$ and
$$g(r)\to-(N-2)\mathfrak c_{k,N} \ \ \ \mbox{as}\ \ r \to0^{+} \ \ \ \mbox{and} \ \ g(r)\to+\infty \ \ \mbox{as} \ \ r\to1^{-}.$$
Consequently:
\begin{enumerate}[label=\textup{(\roman*)}]
\item if $\mathfrak c_{k,N}<0$, then $\Phi_{k,N}'>0$ on $(0,1)$; in particular $\Phi_{k,N}$ has no critical points in $(0,1)$;
\item if $\mathfrak c_{k,N}>0$, then $\Phi_{k,N}$ has exactly one critical point $r_0\in(0,1)$, which is its unique global minimum point, and $\Phi_{k,N}''(r_0)>0$.
\end{enumerate}
\end{proposition}

\begin{proof}
By Lemmas~\ref{lem:decomp} and \ref{lem:series}, differentiating the power series term by term,
\[
g(r)=\sum_{m=1}^\infty 2m\,a_m\,r^{2m+N-2}-(N-2)\,\mathfrak c_{k,N}.
\]
Since $a_m\geqslant0$ and $a_2>0$, it follows that $g$ is strictly increasing and \\
$g(r)\to-(N-2)\mathfrak c_{k,N}$ as $r\to0^+$. \ \\
Moreover, by monotone convergence, 
$$g(r)\to\sum_{m=1}^{+\infty} 2m\,a_m-(N-2)\mathfrak c_{k,N}$$
as $r\to1^-$ and, by Lemma~\ref{lem:R}, we obtain that 
$$\sum_{m=1}^{+\infty} 2m\,a_m\geqslant\sum_{m=1}^{+\infty}a_m=\lim_{r\to1^-}R_{k,N}(r)=+\infty.$$
Moreover, \ \\
$(i)$ if $\mathfrak c_{k,N}<0$, then $g>-(N-2)\mathfrak c_{k,N}>0$ on $(0,1)$; \ \\
$(ii)$ if $\mathfrak c_{k,N}>0$, then $g$ changes sign exactly once, at some $r_0\in(0,1)$, so $\Phi_{k,N}$ is strictly decreasing on $(0,r_0)$ and strictly increasing on $(r_0,1)$. Finally, since $\Phi_{k,N}'(r_0)=0$,
\[r_0^{N-1}\,\Phi_{k,N}''(r_0)=g'(r_0)\geqslant 4a_2(N+2)\,r_0^{N+1}>0.\]
\end{proof}

\begin{corollary}\label{cor:nondeg}
Let $N\geqslant5$ and $2\leqslant k\leqslant k^*(N)$. Then $\Gamma$ has exactly one critical point in $\mathcal M=(0,1)\times(0,\infty)$, namely $(r_0,cd(r_0))$, where $r_0$ is the unique global minimum point of $\Phi_{k,N}$, and it is non-degenerate. In particular, the solutions of Theorem~\ref{main}\,(b) satisfy $r_\e\to r_0$ and $d_\e\to d_0$ along the whole family $\e\to0$.
\end{corollary}

\begin{proof}
It follows from Proposition~\ref{prop:sign}, Proposition~\ref{prop:shape}\,(ii) and Lemma~\ref{lem:critPsi}. The convergence of the whole family follows from the uniqueness of the global minimum point.
\end{proof}

\subsection*{Non-existence for $k>k^*(N)$}

\begin{theorem}\label{thm:nonex}
Let $N\geqslant5$ and $k>k^*(N)$, and let $\mathcal K\subset\mathcal M=(0,1)\times(0,\infty)$ be compact. Then there exists $\e_{\mathcal K}>0$ such that, for every $\e\in(0,\e_{\mathcal K})$, problem \eqref{pb} in $\mathbb{B}_1(0)$ has no solution of the form
\[
u=V_y+\phi,\qquad y\in\mathcal K,\quad \phi\in Z^\perp,\quad \|\phi\|_{H^1_0}\leqslant C\vartheta(\mu),
\]
where $V_y$ is the ansatz \eqref{ansatz} in case \textup{(B)} and $C$ is the constant of Proposition~\ref{prop:reduction}.
\end{theorem}

\begin{proof}
Choose $\eta>0$ such that $\mathcal K\subset\operatorname{int}\mathcal P_\eta$. Assume by contradiction that there are $\e_n\to0$, $y_n\in\mathcal K$ and solutions $u_n=V_{y_n}+\phi_n$ of \eqref{pb} with $\e=\e_n$, $\phi_n\in Z^\perp$ and $\|\phi_n\|\leqslant C\vartheta(\mu_n)$. Since $u_n$ solves \eqref{pb}, the projected equation of Proposition~\ref{prop:reduction} holds, so that $\phi_n=\phi_{y_n}$ by uniqueness. Hence $\nabla\widetilde J_{\e_n}(y_n)=0$, because $\partial_{y_i}\widetilde J_\e(y)=J_\e'(V_y+\phi_y)\bigl[\partial_{y_i}(V_y+\phi_y)\bigr]$. By Proposition~\ref{prop:expansion},
\[
0=\mu_n^{-(N-2)}\nabla\widetilde J_{\e_n}(y_n)=\ell\,\nabla\Gamma_{A,B}(y_n)+\nabla\varpi_{\e_n}(y_n),\qquad \nabla\varpi_{\e_n}\to0\ \text{uniformly on }\mathcal P_\eta .
\]
Up to a subsequence, $y_n\to y^*\in\mathcal K$, and $y^*$ is a critical point of $\Gamma_{A,B}$, hence of $\Gamma$. By Lemma~\ref{lem:critPsi}, $\Phi_{k,N}'(r^*)=0$, which contradicts Proposition~\ref{prop:shape}\,(i), since $\mathfrak c_{k,N}<0$ by Proposition~\ref{prop:sign}.
\end{proof}

\begin{proof}[Proof of Theorem~\ref{main2}]
Part (i) is Corollary~\ref{cor:nondeg}, and part (ii) is Theorem~\ref{thm:nonex}.
\end{proof}

\begin{remark}
Lemma~\ref{lem:series} and Proposition~\ref{prop:shape} hold verbatim for $N=4$, i.e.\ $\lambda=1$. In particular $\Phi_{k,4}$ is strictly increasing on $(0,1)$ for every $k\geqslant7$, since $\mathfrak c_{k,4}<0$ by Remark~\ref{rmk4D}.
\end{remark}

\section*{Data Availability Statements}
All data generated or analysed during this study are included in this article.
\section*{Declarations}
{\bf Conflicts of Interest:} The author declares he has no financial interests.

\end{document}